\documentclass{article}
\usepackage{graphicx}
\usepackage{subfigure}
\usepackage{amsmath}
\usepackage{amssymb}
\usepackage{amsthm}
\usepackage{float}
\usepackage[top=2cm, bottom=2cm, left=2cm, right=2cm]{geometry}

\newtheorem{theorem}{Theorem}
\newtheorem{defn}{Definition}

\usepackage{authblk}

\usepackage[compress]{cite}
\usepackage{xcolor}
\usepackage{enumitem}

\begin{document}

	\title{Deterministic chaos for Markov  chains}	
	
	\author[]{Marat Akhmet\thanks{Corresponding Author, Tel.: +90 312 210 5355, Fax: +90 312 210 2972, E-mail: marat@metu.edu.tr} }
	\affil[]{Department of Mathematics, Middle East Technical University, 06800 Ankara, Turkey}

	\date{}
	\maketitle

	\noindent \textbf{Abstract.}  We   find that   Markov  chains   with  finite  state  space  are  Poincar\'{e}  chaotic.  Moreover,  finite realizations  of the  chains   are   arcs  of  each  unpredictable  orbit   for  sure.  An  illustrating  example   with   a  proper numerical  simulation is provided.

	\section{Introduction} 
	\label{sec1}

	The   main   task   of the  present  research  is  to  find  signs   of    deterministic   chaos in   Markov chains.  Then,   to    show that  the  presence is maximal   in the  sense of probability.  The  method  of domain structured chaos,  \cite{AkhmetDomainStruct},  is in the  basis of our  study. It  has been    successfully  applied for  chaos  indication in fractals \cite{AkhmetSimilarity},  neural  networks \cite{AkhmetDomainStruct},  and for  Bernoulli process \cite{AkhmetBernoulli}.   This time,  we  provide arguments for   the  chaos in  Markov chains.  At  first,    we   formalize the  realizations of the  chains  as  conveniently  constructed sequences.  Then,   give explanations, why the  random dynamics   admits the  chaos features. Considering,  first  of all  the  Poincar\'{e}  chaos. 

	  The  first  goal  of  the  Markov   research  \cite{Markov}   was to  show that  the  random processes of  dependent  events  may  behave as   processes with independent events.    Thus,   the    model,    simple   and most  effective for  many  applications  was invented.   It  is   impossible underestimate the  role of  the Markov chains and processes   in development  of random dynamics  theory  and its applications.  For  instance,  the egrodic  theorem was  strictly  approved   at  the first time under the  circumstances.    There are  many  observations that  the  chains  are  relative  to  symbolic  dynamics and correspondingly to  Bernoulli   scheme.   The significant   step  in  the  comprehension    was done   by   Donald Ornstein in  \cite{Ornstein1},   who  verified that   $B-$automorphisms  such  as sub-shifts of finite type and Markov shifts, Anosov flows and Sinai's billiards, ergodic automorphisms of the $n-$torus  \cite{Ponce}, and the continued fraction transform    are,  in fact,  isomorphic.   Issuing    from these  sources  it  is  of   great  interest  not  only  to show that  different    stochastic processes  can be described in the terms of deterministic  chaos,  but   that  they  relate   equally  to   the  chaos.   Researchers  have looked  in both  directions,   for  chaos  in random dynamics, as well  as for stochastic features in deterministic  motions \cite{Alekseev,Bowen,Bowen1,Bowen2,Bowen3,Bowen4,Bowen5,Ornstein2,Ornstein3}.   Thus,  the problem  of  chaos  in Markov chains,  which  is discussed in the paper,    is a part  of the  more  general  project.  	One can apply   finite \cite{Alekseev,Bowen} or countable  \cite{Bunimovich}  partitioning to  find  randomness in  deterministic  motions.    Differently,  in the   present research   we utilize an   uncountable  partitioning  to   establish   strong  relations of the  random processes   with   deterministic   chaos.

	   In   resent  studies  \cite{AkhmetSimilarity,AkhmetDomainStruct,AkhmetAbstFract,AkhmetCube, AkhmetBook}, we formalize   a   chaos  generation through  specially   structured   sets,   subduing  them   to  the  similarity  map,  which  is  in the  paradigm  of the  Bernoulli  shift \cite{Wiggins88}.   This   provides  to  us  the  strong  comprehension   that  the chaotic behavior   is proper  to  the  Bernoulli scheme   dynamics and then  to   various  random processes on discrete    and continuous time \cite{AkhmetHyperbolic}-\cite{AkhmetModular}.
	    Moreover,  it  was clarified   in which  sense one should   understand presence of deterministic chaos  attributes in stochastic  processes.   
	    
	   In the  present  research  we are  on the  next  step  of the   exploration,  and show that  the chaos  can  be recognized for  Markov chains.    Time-homogeneous  chains   on  a  finite state space,  and chains with  memory  are discussed.   The  approach  is  applied  such  that   the dynamics  can  be investigated,  focusing  not  on ergodicity,  but  on a single   motion description such  that    an  individual   geometry  is  better seen.    Thus,  we  show  that    an isolated realization of a Markov chain behaves in time identically  with  a properly  chosen  realization of a Bernoulli  scheme,  and  they both  are identical  in dynamics with   the  path  of   the  similarity  map  in  a  correspondingly  chosen   space.  The  unpredictable  orbit \cite{AkhmetUnpredictable}  as a single  isolated motion, presenting  the  Poincar\'{e} chaos  \cite{AkhmetPoincare}  in the chains,   is  the  ultimate  point   of   the  comprehension.   
	   
	   Considering  the  signs  of the  chaos  in Markov chains,   it  is important   to  say  not  about  sensitivity,  but  individual   behavior  of trajectories.   That  we say  that  simulating   the  random process,   one can  not  determine that  it  is different  than  an  orbit  of corresponding    deterministic   dynamics,  which  we can  suggest   for  sure.  But  discussing  the  sensitivity,  one can  say  that   it  is true with  some positive probability.  More precisely,    simulating   a random process,  we are confident that  there  is   another  motion   which  starts 
	   arbitrary  near   and  diverges from the first  one   on  the  distance not  less than  a positive  constant   common for  the  dynamics with  a positive probability.  Of course,   the  probability   less than  one can  make   us critical,  but  it  does not  different  than that  for  the  deterministic chaos  as  the   points , which  are  start for  the  diverging can  not  be determined for  both  dynamics, deterministic as well  as for the  random.  But,  we in our research  know how to  find  the   start  point.  For this we have developed the  sequential  test. On this  basis, one can not  only  evaluate  the  sequences, which are in the  basis  of the  definition of the  unpredictable  point, but   also determine point   which   provide the  divergence with   arbitrary  smallness.  This  everything  is done on the  basis of a single orbit  or  let  say  a single simulation.

	   \section{Preliminaries}
	   
	   The  domain structured  chaos unites models,  which  admit   all  types of theoretical  chaos such  as Poincar\'{e},  Li-Yorke and Devaney \cite{AkhmetDomainStruct},  but  the  first  one is most  convenient  to   confirm   that   chaos  is  proper for  the  random dynamics.   For this  reason, let  us start  with  the  definition of the  unpredictable  sequence.   The  notion is  the  strict  evidence  for  the  Poincar\'{e}  chaos \cite{AkhmetPoincare}.

	   \begin{defn} \cite{AkhmetPoincare} \label{defseq}
	   	A bounded sequence  $ \kappa_i, i  \ge 1,$ 	in $ \mathbb R^p $ is called unpredictable if there exist a positive number $ \varepsilon_{0} $ and sequences $  \zeta_n $, $ \eta_n $, $ n   \ge 1$, of positive integers both of which diverge to infinity such that $ \| \kappa_{i + \zeta_n} - \kappa_i \| \to 0 $ as $ n \to \infty$ for each $ i $ in bounded intervals of integers and $ \| \kappa_{\zeta_n + \eta_n} - \kappa_{\eta_n} \| \geq \varepsilon_{0} $ for each   natural  number  $ n.$
	   \end{defn}

	         For  a  finite  metric  space $(S,d),$  the  last  definition   has   the  next  form.

	   \begin{defn} \label{defseq1}
	   	A bounded sequence $ s_i\in S,  i  \ge 0,$ is called unpredictable if there exist a positive number $ \varepsilon_{0} $ and sequences $  \zeta_n $, $ \eta_n ,  n \ge 0,$ of positive integers both of which diverge to infinity such that  $s_{i + \zeta_n} =  s_i$   for each   bounded interval of integers,    if $n$   is sufficiently  large,   and $ d(s_{\zeta_n + \eta_n}, s_{\eta_n})\geq \varepsilon_{0} $ for each   natural  number $ n.$
	   \end{defn}

	   The    main result   of the  paper is  that  there   exists  the   realization,  an unpredictable    sequence,   of  the  Markov process  with  the  finite  state  space,  which  closure   in the  topology   of convergence on bounded intervals  is the set  of all  infinite realizations.    Then,  it  implies  that   each   finite   simulation of the  dynamics is an ark  of  the unpredictable  realization.  This  is why, one can  confirm that \textit{ the Poincar\'{e} chaos is   a certain   event    for   the  Markov chain}.   Let  us remind that  an event is   certain, if   it  occurs  at every performance of an experiment.  One   must  specify  that  the  convergence,  in the circumstances,  is the  coincidence  of the  unpredictable  realization with  each  finite  realization.

	   The  closure  of the unpredictable  realization is said to  be the  quasi-minimal  set \cite{Sell}. It  contains uncountable  set  of unpredictable  realizations.    We have proved    that  \textit{a quasi-minimal  set  as the  union of all infinite  realizations of the  Markov  chain is  a   certain  event}. And this  is another  formulation of the main result.  The results  on  chaos  presence in the  random processes  have been  investigated  in our  previous  papers 
	   \cite{AkhmetBernoulli,AkhmetRandom},  but  the   the chaos  appearance  in the  dynamics with  probability  one is  approved   at   the  first  time.     We have explored  that  namely unpredictable  orbit  is  most  proper   candidate  for  the  analysis.

	   Generally speaking,    our  results  confirm that  the  Bernoulli  scheme \cite{AkhmetBernoulli},  Markov chains  as well  as abstract  hyperbolic  dynamics \cite{AkhmetHyperbolic}  are  all   with  the  same type of chaos.  It  is    significant   that  the  Poincar\'{e}   chaos  is proper for  the motions.  This  provides,  new opportunities,   exceptionally  for    stochastic  processes.   We    suppose   that  the   research  can  be complemented  with    similar  analysis  for other  Bernoulli  automorphisms   considered  by  D. Ornstein  \cite{Ornstein1} as well  with  extension of the  results  for majority  of stochastic  processes,  if  proper structured domains will be  constructed.

	A Markov chain is a stochastic model,   which   describes a sequence of possible events   such  that   the probability of each event depends only on the state attained in the previous   one \cite{Karlin,Hayek}.  There  are   many applications  of the Markov chains as statistical models of real-world processes such as studying  queues or lines of customers arriving at an airport, currency exchange rates,  cruise control systems in motor vehicles   and animal population dynamics \cite{Meyn}.
	
	Consider   a  discrete-time stochastic process  $ X_n , n \ge 0,$ on a  countable   set 	$S.$    That  is,  a collection of  random variables defined on a probability space
$(\Omega,F, P),$  where  $ P$   is a probability measure on a family of events  $F$ 
in an event-space $\Omega.$  The set $S$  is the state space of the process, and the
value $X_n  \in  S$  is the state of the process at time $n.$ 

The   Markov chain,   is   a   stochastic  process  such  that  
   the   Markov property    $P\{X_{n+1} = s_j|X_0, . . . , X_n\} = P\{X_{n+1} = s_j|X_n\}$   is  true  for   all $ s_i,s_ j \in  S$ and $n \ge 0,$    and   $P\{X_{n+1} = s_j|X_n = s_i\} = p_{ij},$ where   $p_{ij}$  is the   transition  probability  that  the  chain jumps from state  $i$  to  state $j.$ 
The property says that, at any time $n,$  the
next state $X_{n+1}$  is conditionally independent of the past $X_0, . . . , X_{n-1}$  given
the present state $X_n.$     More precisely,  that  the transition probabilities do not depend on
the time parameter $n.$   That  is,   the chain is time-homogeneous. If
the transition probabilities were functions of time, the process  would be a
non-time-homogeneous Markov chain. 

  To    avoid any   terminological  confuse,  we  will appeal  to  the  dynamical interpretation  of   Markov chains  as   the  following  recursion
  
  \begin{equation}  \label{dynamics}
  X_n = f(X_{n-1},Y_n),
  \end{equation}  
   where   $n \ge 0$    is the  time parameter,   $Y_1,Y_2 ... $ are independent and identically distributed  and $f$  is a deterministic function. That is, the new state $ X_n $ is simply a function of the last state and an auxiliary random variable.     In  other  words,  one   can  consider  a Markov chain  as a random   dynamical  system \cite{Arnold}.  In the  present paper    realizations   are   considered orbits   or  trajectories of  the     corresponding random  dynamics   (\ref{dynamics}).   Consequently,    one  can  denote  the  realization as $X_n = X(n) = X(n,X_0),$   where  $X_0 = X(0)$  is the  initial  value,  which   is determined randomly.    This  is   a step  of  better  comprehension of the  stochastic  dynamics  through    chaotic  interpretation.  If  the   paparmeter  $n$ runs over $\mathbb N$  we say     that  a   realization  is  infinite.   Otherwise,  it s  finite.    Thus,  realizations are  infinite  or finite  sequences of elements from the  state  space.   We   will  use the  set  of all  realizations to  discuss the  problem  of the   deterministc   chaos  for  the  Markov chains. 
   

\section{Domain structured  chaos and Markov chains} 
	
		Consider a finite state   space $S= \{s_1,\ldots,s_m\},$  where  $m$  is   a natural  number,   not  smaller than  two,     and    a  metric $d$  for  the  space.    		Denote  $p_{ij} = p_j(s_i)$   the  Markov probability    for  $s_j, j = 1,\ldots,m,$  such  that     $\Sigma_{j=1}^m p_{ij} =1$    for  all  $i =1,\ldots,m.$  Assume  that   all     probabilities $p_{ij}, i,j = 1, \ldots,m, $  are  positive.
		
In the  basis  of our  construction is  the  event  $f_{ij},  i,j= 1,\ldots,m,$   which  consists  of two     elementary  events,     $s_i$  and $s_j,$    happen successively,   such  that     an  infinite  realization of the    Markov chain   is   formalized as   the   infinite   sequence $f_{i_1j_1}f_{i_2j_2}\ldots f_{i_nj_n}\ldots,$  with  $j_k = i_{k-1}$   for  all 	$k =2,3,\ldots.$      We   have that  $p(f_{ij})= p_{ij}.$    The   formalization does not  give advantages,  if one consider the  chains without  memory,  but  it   makes  easier  the  discussion of the  processes with  memory,  in what  follows.  Present  the last  sequence as the  element   $\mathcal{F}_{i_1 i_2 \ldots}, i_k=1,2, ..., m, \; k=1, 2,\ldots,$  of the      space $\mathcal{F}$ with  metric  $\delta(\mathcal{F}_{i_1 i_2 \ldots},\mathcal{F}_{l_1 l_2 \ldots}) = \Sigma_{k=1}^{\infty} d(s_{i_k},s_{l_k})/2^k.$    
			
		Next,   we   formalize  Markov chains   with    memory of  a non-zero  length .  We start  with   the  length  equal  to  two such  that  the  element  $f_{ijk}$  presents   three elementary  events  $s_i,  s_j$  and $s_k$ happen   successfully,  and   the  probability      for  $s_k$  is equal  to  $p_{ijk} = p_k(s_i,s_j).$   Then the  Markov chain  with  memory    has the  formal  presentation    $f_{i_1j_1k_1}f_{i_2j_2k_3}\ldots f_{i_nj_nk_n}\ldots,$  where   $j_l = k_{l-1},i_l=j_{l-1}$   for  all 	$l =2,3,\ldots.$    Accepting   the  last  sequence  as   the  element   $\mathcal{F}_{i_1 i_2 \ldots},  i_k=1,2, ..., m, \; k=1, 2,\ldots,$  of the      space $\mathcal{F}$ with  metric  $\delta(\mathcal{F}_{i_1 i_2 \ldots},\mathcal{F}_{l_1 l_2 \ldots}) = \Sigma_{k=1}^{\infty} d(s_{i_k},s_{l_k})/2^k$    we  attain the  basis,  common   with  that  for    the  chain  without  memory. 
		  
		  At  last,  consider the Markov process  with   the  memory  of  the   length equal  to  arbitrary  natural  number $n.$  Then   we formalize the  discussion with    the    elements  $f_{i^1,\ldots,i^n}, $  which  consist of   successive elementary  events  $s_{i^1},\ldots,s_{i^n},$  such  that   the  sequence $f_{i_1^1,\ldots,i_1^n}\ldots,f_{i_l^1,\ldots,i_l^n}\ldots$  with  $i_l^j = i_{l-1}^{j+1},   j=1,\ldots,n-1, l =2,3,\ldots,$  is the  formalization of the  chain.    We  obtain  the  structure    for  the  dynamics  research,  if    accept  the  last  sequence    as an  element  		  
		   $\mathcal{F}_{i_1^1i_2^1 \ldots},  i_k^1=1,2, ..., m, \; k=1, 2,\ldots,$   of the  space  $\mathcal{F},$  making  stress on the events with  the indices   $i_k^1, k =1,2,\ldots.$        To  complete  the  chaos  analysis  we shall  consider,   other  indexes  also,    namely,  the   spaces   consisting  of  elements   $\mathcal{F}_{i_1^ji_2^j \ldots},  i_k^j=1,2, ..., m, \; k=1, 2,\ldots,$    with  arbitrary  fixed   $j =1,2,\ldots,n.$    This  is  not  necessary,   if one   consider the    dynamics   from the  traditional  point  of view,  when   the  phenomenon has to  be observed only  for  unbounded sequence of   moments, even  for continuous  time. Remember, Poincar\'{e}    stroboscopic approach.    Nevertheless,  in the  research  of the  Markov chains   with  memory,    that  is  in our  present  case,  it  is  significant  to  precise that   there  are   finite  number of   subsequences which  cover   with  chaotic  dynamics   over the  whole  discrete  time range.    Thus  the  specific properties of  the  dynamics  are  emphasized.
		   
		    It  is  clear  that  for  all  cases,  regardless are  they   with  memory  or   without  memory,   we have constructed one and the  same    space  $\mathcal F$ of elements  $\mathcal{F}_{i_1i_2 \ldots},  i_k=1,2, ..., m, \; k=1, 2,\ldots,$  with  the  distance
		   $\delta(\mathcal{F}_{i_1 i_2 \ldots},\mathcal{F}_{l_1 l_2 \ldots}) = \Sigma_{k=1}^{\infty} d(s_{i_k},s_{l_k})/2^k.$      This   is why,  the  space   is    the   object  of analysis    for  the  deterministic  chaos  presence,  next.    To   complete  the dynamics,    we shall  need the  special  map  on the space.  
		   
		   	 Consider the map $ \varphi : \mathcal{F} \to \mathcal{F} $ such that
			\begin{equation} \label{MapDefn}
			\varphi (\mathcal{F}_{i_1 i_2 ... i_n ... }) = \mathcal{F}_{i_2 i_3 ... i_n ... },
			\end{equation}		
			for  each  element  of the  set.  The map $\varphi$  is   in the  paradigm  of  the   Bernoulli shift \cite{Wiggins88},   known    for the symbolic dynamics. It  is said to  be  the  \textit{similarity  map} \cite{AkhmetSimilarity} as it is   convenient   to  describe  fractals,  which  are determined through  the  self-similarity. 
			
			From the  definitions of  the  set $\mathcal F$  and  map  $\varphi$   it  implies that  the  values  of  the map      correspond   to  the    members of  the  state  space,  which  appear  orderly.  Consequently, if one   proves that   the   process  of appearance  is chaotic  in one of the   senses accepted in literature,  then  we   must  accept  that the  chaos  is proper for  the  stochastic dynamics  with  proper   arguments  of  probability.

		The following  sets are  needed, 
		\begin{equation} \label{AbstFracSubSet}
		\mathcal{F}_{i_1 i_2 ... i_n} = \bigcup_{j_k=1,2, ..., m } \mathcal{F}_{i_1 i_2 ... i_n j_1 j_2 ... },
		\end{equation}
		where indices  $ i_1, i_2, ..., i_n,$  are fixed.
		
		   It is clear that
		   \[ \mathcal{F} \supseteq \mathcal{F}_{i_1} \supseteq \mathcal{F}_{i_1 i_2} \supseteq ... \supseteq \mathcal{F}_{i_1 i_2 ... i_n} \supseteq \mathcal{F}_{i_1 i_2 ... i_n i_{n+1}} ... , \; i_k=1, 2, ... , m, \; k=1, 2, ... \, ,\]
		   that is, the sets form a nested sequence. 
		   
		   	Considering iterations of the map, one can verify that
		   \begin{equation} \label{MapSubset}
		   \varphi^n(\mathcal F_{i_1 i_2 ... i_n}) = \mathcal{F},
		   \end{equation}
		   for arbitrary natural number $ n $ and $ i_k=1,2, ..., m, \; k=1, 2, ... \, $. The relations (\ref{MapDefn}) and (\ref{MapSubset}) give us a reason to call $ \varphi $ a \textit{similarity map} and the number $ n $ the \textit{order of similarity}.
		   
		   We will  say   that for the sets $ \mathcal{F}_{i_1 i_2 ... i_n} $ the \textit{diameter condition} is valid,  if 
		   \begin{equation} \label{Diamprop}
		   \max_{i_k=1,2, ..., m} \mathrm{diam}(\mathcal{F}_{i_1 i_2 ... i_n}) \to 0 \;\; \text{as} \;\; n \to \infty,
		   \end{equation}
		   where $ \mathrm{diam}(A) = \sup \{ d(\textbf{x}, \textbf{y}) : \textbf{x}, \textbf{y} \in A \} $, for a set $ A $ in $ \mathcal{F} $.

		    Denote the distance between two nonempty bounded sets $ A $ and $ B $ in $ \mathcal{F} $ by $ d(A, B)= \inf \{ d(\textbf{x}, \textbf{y}) : \textbf{x}\in A, \, \textbf{y} \in B \} $. Set $ \mathcal{F} $ satisfies the \textit{separation condition }of degree $ n $ if there exist a positive number $ \varepsilon_0 $ and a natural number $ n $  such that for arbitrary   indices  $ i_1 i_2 ... i_n $ one can find   indices $ j_1 j_2 ... j_n $ such  that
		   \begin{equation} \label{C2}
		   d \big( \mathcal{F}_{i_1 i_2 ... i_n} \, , \, \mathcal{F}_{j_1 j_2 ... j_n} \big) \geq \varepsilon_0.
		   \end{equation}

     		   Next,   we will  formulate  two  theorems on   Poincar\'{e}  and Devanye   chaos. Verification of the  assertion is given in \cite{AkhmetBernoulli}. 
     		    The following theorem  asserts  that   the similarity map $ \varphi $ possesses the three ingredients of Devaney chaos, namely density of periodic points, transitivity and sensitivity. A point $ \mathcal{F}_{i_1 i_2 i_3 ...} \in \mathcal{F} $ is periodic with period $ n $ if its index consists of endless repetitions of a block of $ n $ terms.
		   If    t   said to be  a \textit{chaotic  structure}  for   $F.$  The  next  theorem    is  a particular  case of more  general  Theorem  \ref{Thm3}.
		   
		   \begin{theorem} \label{Thm1}  \cite{AkhmetBernoulli, AkhmetDomainStruct} If   the   diameter and separation conditions are valid,  then  the  dynamics $(\mathcal F, \delta, \phi)$  is chaotic in the sense of Devaney.
		   \end{theorem}

		   In \cite{AkhmetUnpredictable,AkhmetPoincare}, Poisson stable motion is utilized to distinguish  chaotic behavior from  periodic motions in Devaney and Li-Yorke types.  The dynamics  is given the named  Poincar\`{e} chaos.  The next theorem shows that the Poincar\`{e} chaos is valid for the similarity dynamics.
		   
		   \begin{theorem} \label{Thm2}   \cite{AkhmetBernoulli,AkhmetDomainStruct}  If   the   diameter and separation conditions are valid,   then the  similarity map possesses Poincar\`{e} chaos.
		   \end{theorem}

		   In addition to the Devaney and Poincar\`{e} chaos, it can be shown that the Li-Yorke chaos is also present  in the dynamics of the map $ \varphi $. The proof of the  theorem is similar to that of Theorem 6.35 in \cite{Chen} for the shift map defined in  the space of symbolic sequences.

		   Thus,  it is proven   that   the dynamics admits the ingredients  of all  the  theoretical  chaos,  and  we say  that     		   
		 that  for  the  space $(\mathcal F, \delta, \varphi)$ \textit{the  domain structured chaos}  is proper.    For the  present  research  it  is  important   that  there   exists  an   unpredictable trajectory   of $\varphi$   in the  sense of Deefinition \ref{defseq1}.     It  is Poisson stable,  and the   closure   of the  set  of all  trajectories is a quasi-minimal set.    In our  research  \cite{AkhmetUnpredictable}  we have proved that  if a Poisson stable point   is additionally  unpredictable,  then there  is the  sensitivity  in the  set  of motions.     Thus,  it  was recognized that  there  is Poincar\'{e}  chaos.  The set  $\mathcal{F}$  is bounded.  The  convergence is in the  topology  on bounded sets.   The  orbits  of the  map  $\varphi$  are  infinite  sequences.
		 
		   Since of the  accordance  between  the  Markov chain on the  finite  state   space and the  set  dynamics of  the  map  $\varphi,$  one can conclude that  the  following  assertion is valid.  To  formulate  the  result, let  us fix  an  unpredictable  realization of the  chain, which  can  be determined as follows. 
		   Consider an unpredictable  point  $\mathcal F^*_{i_1i_2\ldots}$  of the map  $\varphi.$  Fix the  sequence,  $s^* = \{ s^*_{i_k}\}_k,$   which  is the corresponding realization.  It  is an unpredictable sequence.      Duo  to  the  Definition   \ref{defseq1}   and  Theorem  \ref{Thm2},  the  following  assertion  is valid.
		   
		   \begin{theorem}  \label{Thm3} Each  finite  realization of the  Markov chain  coincides with  an   arc of    $s^*.$   That  is,  the realization   happens   in each  experiment of the   chain,  and  is a certain event. 
		   \end{theorem}
	   
	 The   sensitivity   property  is obvious,  since  of the  finite  state  space to  start  arbitrary  near  means,  to coincide    at  the start  moment, and   absense   of   the  divergence   means absence of the  randomness.

		      Thus,  one has to  recognize    that  each  experiment with  the  Markov chain   (without  memory) produces an arc of an  unpredictable orbit.    This  result  can  be considered  as the main one in this reseach. It   is in full accordance   with  the  principle of the  ergodic  theory  \cite{Walters}  that  a single trajectory  proves behavior  of the  whole dynamics and all  other trajectories.  In fact,  we can  say  that  it  is a fixed unpredictable  orbit.   In other words,  this is reproduction of the  chaos  with  probability  one.    This   is what   guaranties the  irregularity  of each  finite  sample  path  of the  chain.  As it  was proved \cite{AkhmetPoincare}, the  existence   of   an  unpredictable  orbit  implies sensitivity. Consequently,   for  each  infinite  sample path  one can  find arbitrary  near   another  realization, which  definitely  diverges at  finite  time    for some positive  number common in the  dynamics.   Additionally,   our  research  confirms that  constructing  a    finite  realization of the  chain numerically,  we build a piece  of the  graph  of unpredictable   sequence,  and this  can  be applied for  definition  of  unpredictable  functions. This  was realized in our  paper \cite{AkhmetRandom}.
		      Finally,    consider     a fixed infinite   realization.   One  can see that    arbitrary  near  to  the  sequence    a finite realization  starts,  which  exists  with  positive probability and diverges   from the    fixed  realization on the  distance not  less than  a positive constant  common for  the chain. That  is the  sensitivity   presents  with  non-zero  probability. We  suppose that  the  probability is  equal  to  the  unit,  but  this  is  problem  for  the  future  research. 
		      
		      Now,  let  us  focus  on the  Markov chains with  memory.   Let  the  space  $\mathcal{F}$  of       
		       elements  $\mathcal{F}_{i_1^1i_2^1 \ldots},  i_k^1=1,2, ..., m, \; k=1, 2,\ldots$  is given    with  the  similarity  map  $\varphi.$
		       It  is clear that the  Poincar\'{e}   chaos   with  probability one   is  present in the  discrete  dynamics.   The  discussion  is identical  with  that    done for  the  chains without  memory.   Consequently, one can  conclude that  the  chaos exists  in the  dynamics   with probability  one.  If  one objects  the  assertions,   since  we made the decision considering   subsequences of realizations, our  response  is that  this  is true  for  many  chaos research  in literature. For  example,  we  indicate   chaos     for  continuous  dynamics just  by  considering Poincar\'{e}  sections observation.  In our  case, the  arguments are  much  more  strong,  as we observe the  chaos  for  all sets    $\mathcal{F}_{i_j^1i_2^j \ldots},  i_k^1=1,2, ..., m, \; k=1, 2,\ldots, j =1,\ldots,m.$   Consequently, 
		       we can  make  decision that  the  chaos  is more  strong  in its  presence  for  the  chains with  memory.

		      	 \subsection{An example: random walk} 
		      	 
\begin{figure}[ht]
		      	   	      	 	\centering 
		      	   	      	 	 \includegraphics[width=0.6\textwidth]{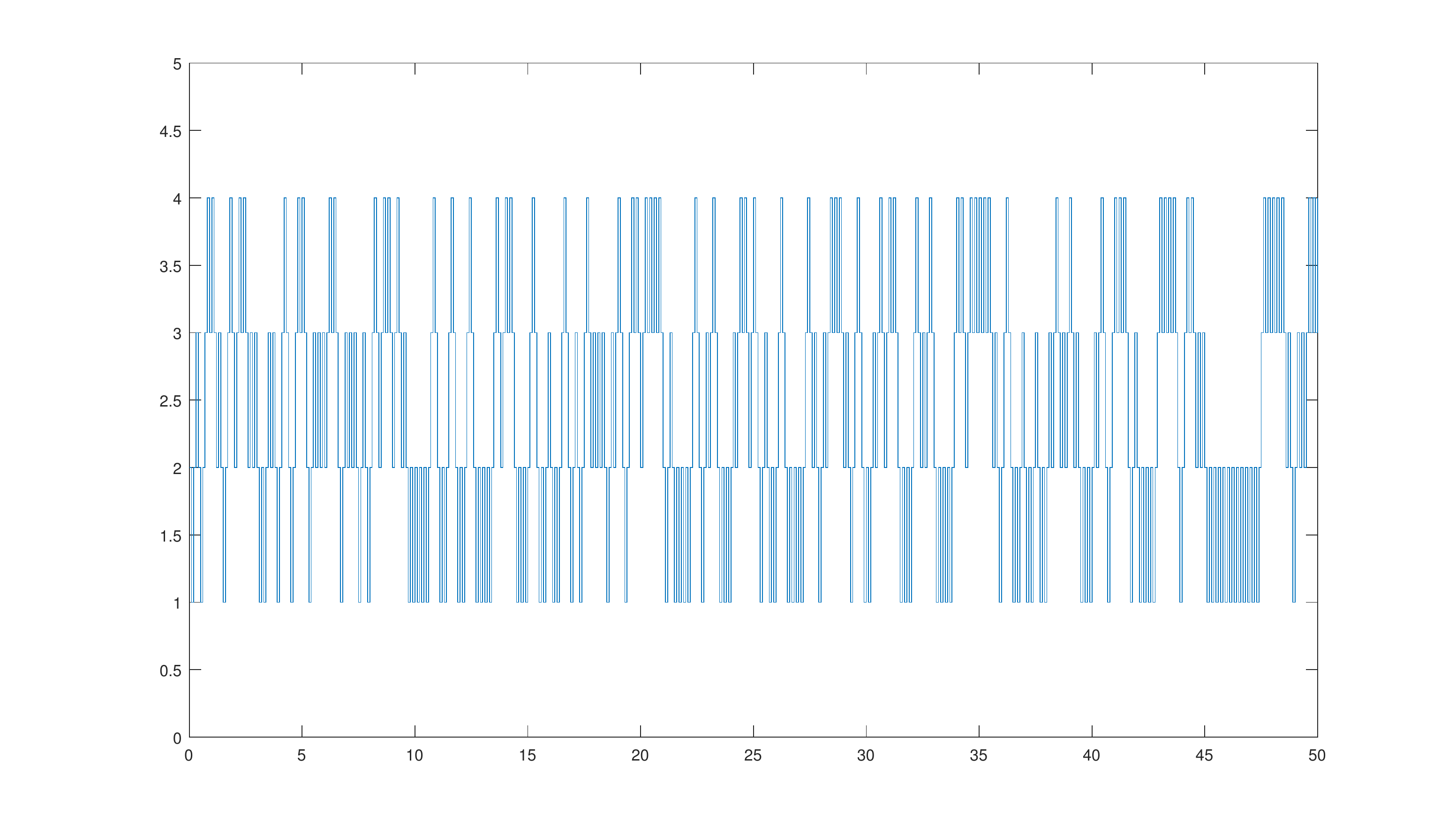}
		      	   	      	 	\caption{The  graph  of the  function $\phi(t)$  as visualiation of the   unpredicatble  realizatoin   of the  chain $X_n, n \ge 0.$  The  vertical  lines   connecting  the  pieces of the  graph are drawn,  for better  visibility.}
		      	   	      	 		      	 	\label{f1}
		      	   	      	 \end{figure}

		      	  Consider,  as  an   example, the  following  Markovian  chain.   Let  the    real  valued scalar  dynamics $X_{n+1}= X_n + Y_n,  n  \ge 0,$  be given  such  that   $Y_n = \{-1,1\}$  is a random variable, with probability  distribution $P(1) =P(-1) = 1/2,$  if $X_n \neq  1,4,$   and $Y_n = -1,$  if $X_n = 4,$   and $Y_n = 1,$ if $X_n = 1.$     To    satisfy  the  construction of the  present research,   we  will  make  the  following   agreements.   First  of all,  denote $s_0 =1, s_1 = 2, s_2= 3, s_3 = 4.$    Consider,  the  state  space $S =\{s_1,s_2\}.$     Introduce   the  following  events,  $f_{11} = \{s_1,s_0,s_1\}, f{12}=\{s_1,s_2\}, f_{21} =\{s_2,s_1\}, f{22} =\{s_2,s_3,s_2\}.$  It  is clear   that   $p_{i1} +p_{i2} =1, i =1,2,$   and all  the  probabilities   are equal  to  the  half.   That  is,     we   are  in the circumstances of the  theory of the  present  paper.     Consequently,  there is  Poincar\'{e}  chaos. To  visualize    an  unpredictable  realization,   we will  draw the  graph  of  the  function $\phi(t) = X_n,$   if  $t \in [n/10,(n+1)/10),   0 \le t\le  60.$   According  to  the  last  Theorem, it  is an  arc of an  unpredictable  sequence.  The  graph of  the  function is  seen  in Figure  \ref{f1}.      It  illustrates   an  unpredictable sequence,  the  sample  path  of the  random  walk.   For  better  visibility of the dynamics   the  vertical  lines   connecting  pieces of the  graph are  drawn.

\section{Conclusion}

The  outcome   of the  research  is  the   existence   of   an  unpredictable   sequence as a realization  of the  Markov  chain, and the    sequence  appears   as  finite  realization   of   each  experiment  of the  process.   That  is,  appearance of the sequence is a certain  event  for  the  stochastic  dynamics.   From this point  of view  one   can  say  that    the  deterministic  chaos  is a certain event  for  the  stochastic dynamics.     This  result  is true  for many  other  discrete  time  random processes. For instance,   the  Bernoulli  scheme.     The  significant  use of the   investigation is that  one can  unite   methods of deterministic  chaos  with  those   for  stochastic dynamics.   Many  other  opportunities may  appear.  Among  the  methods  are  controllability  and  synchronization of chaos
\cite{Gonzales} as well as   different  ways   of chaos  generation \cite{Feckan,Luo}.     We   have  proved  the  sensitivity  is present  certainly.   Thus,  the  deterministic  chaos  has been  approved for  the  stochastic processes.  Evidently,  it  is true for  the  Bernoulli  scheme \cite{AkhmetBernoulli}  and other   dynamics,  wcich   can  be approved for  the  domain structuring. 
Nex our  study  will  relate  with  Markov processes   with  continuous  time,  as well  as   unbounded.  This  also  relates to  many  other  random processes. 
Our  results  provide more lights  on the  Markov chains  as ergodic processes, since   we have shown  that    there  is   the  uncountable  set  of realizations,  unpredictable orbits  and each  of them are  dense in the  set  of all  realizations.

\end{document}